\documentclass[10pt,reqno]{amsart}

\usepackage{titlesec}

\titleformat{\section}
  {\normalfont\fontsize{11}{15}\bfseries}{\thesection}{1em}{}
  
  \titleformat{\subsection}
  {\normalfont\fontsize{10}{15}\itshape}{\thesubsection}{1em}{}

\usepackage{mathrsfs}
\usepackage{amssymb}
\usepackage{amscd}
\usepackage{amsfonts}
\usepackage{version}

\usepackage{graphicx}

\newtheorem{theorem}{Theorem}[section]

\newtheorem{conjecture}[theorem]{Conjecture}

\theoremstyle{definition}
\newtheorem{definition}[theorem]{Definition}

\newtheorem{example}{Example}

\numberwithin{equation}{section}

\def\C{\mathbb{C}}
\def\E{\mathbb{E}}
\def\Z{\mathbb{Z}}
\def\R{\mathbb{R}}
\def\F{\mathbb{F}}
\def\eps{\varepsilon}
\newcommand\id{\operatorname{id}}
\newcommand\SL{\operatorname{SL}}
\newcommand\PSL{\operatorname{PSL}}
\newcommand{\heis}[3]{ \left(\begin{smallmatrix} 1 & #1 & #3 \\ 0 & 1 & #2 \\ 0 & 0 & 1 \end{smallmatrix}\right)  }

\renewcommand{\leq}{\leqslant}
\renewcommand{\geq}{\geqslant}

\author{Ben Green}
\thanks{The author is supported by ERC Starting Grant number 274938 \emph{Approximate algebraic structure and applications.}}
\address{Mathematical Institute\\
Radcliffe Observatory Quarter\\
Woodstock Road\\
Oxford OX2 6GG\\
England }
\email{ben.green@maths.ox.ac.uk}

\title[Approximate algebraic structure]{Approximate algebraic structure}

\begin{document}

\begin{abstract}
We discuss a selection of recent developments in arithmetic combinatorics having to do with ``approximate algebraic structure'' together with some of their applications.
\end{abstract}

\maketitle

\section{Introduction}

Given an inequality, an extremely natural question to ask is
\[ \mbox{When does equality occur?}\]
If a satisfactory answer to this is available, one might then ask
\[ \mbox{When does equality \emph{almost} occur?} \]
To be a little more precise, suppose that we have some family of functions $\mathscr{F}$ and some map (functional) $v : \mathscr{F} \rightarrow \R$. The inequality we are considering might then be of the form
\[ v(f) \leq M \quad \mbox{for all $f \in \mathscr{F}$}.\]
To give an example, the well-known isoperimetric inequality on $\R^n$ may be stated in this form, with $\mathscr{F}$ being the set of all functions $1_A$ where $A \subset \R^n$ is bounded and open (say), $v(1_A)$ being the isoperimetric ratio $\frac{|A|}{|\partial A|^{n/(n-1)}}$, and $M$ being the isoperimetric ratio of any Euclidean ball. 

The first natural question, the \emph{equality question}, is then
\[ \mbox{For which $f \in \mathscr{F}$ do we have $v(f) = M$?}\]
In the case of the isoperimetric inequality, it is well-known (and invariably stated as part of the inequality) that $f = 1_A$ with $A$ being a Euclidean ball. 

The second natural question is 
\[ \mbox{For which $f \in \mathscr{F}$ do we have $v(f) \approx M$?}\]
Of course, to make this precise we must specify what is meant by $\approx$. We further distinguish between what might be called the \emph{stability question}, which asks
\[ \mbox{For which $f \in \mathscr{F}$ do we have $v(f) \geq (1 + o(1)) M$?}\] and what I shall term the \emph{robustness question}, which asks
\[ \mbox{For which $f \in \mathscr{F}$ do we have $v(f) \geq \frac{1}{100}M$ (say)?}\]
Most of this article will be concerned with the robustness question for two particular inequalities, an instance of \emph{Young's inequality for convolutions} and an inequality concerning the \emph{Gowers norms}. In both situations the equality cases are easily established and are highly algebraic in nature (essentially they characterise finite groups and polynomial phases respectively). In both cases study of the robustness question has proven to be surprisingly subtle and has led to diverse applications in areas as different as group theory, additive prime number theory and theoretical computer science. 

The stability question for these same inequalities is much better understood, though it is still nontrivial and has many applications. For want of space, we will not say a great deal about it. As it turns out the stability question for the isoperimetric inequality and related inequalities such as the Brunn-Minkowski inequality is the subject of much current research, not entirely unrelated to the topics discussed in this article: see for example \cite{figalli-jerison1, figalli-jerison2}.

\section{Approximate groups}

\subsection{Young's inequality} Let $G$ be a group with identity element $\id_G$, and let $\mathscr{F}$ be the collection of all finitely-supported functions $f : G \rightarrow [0,\infty)$ with $\sum_{x \in G} f(x) = 1$, $f(x) = f(x^{-1})$ for all $x$ and $f(\id_G) > 0$ . One may think of $f$ as a probability measure on $G$, the measure of a set $A \subset G$ being $\sum_{x \in A} f(x)$.  A particular (rather simple) case of a well-known inequality of Young \cite{young} for convolutions is the bound
\[ v(f) \leq 1 \qquad \mbox{for all $f \in \mathscr{F}$},\]
where
\[ v(f) = \frac{\Vert f \ast f \Vert_2^2}{\Vert f \Vert_2^2} = \frac{\sum_{x \in G} (\sum_{y \in G} f(y) f(yx))^2}{\sum_{x \in G} f(x)^2}.\]
Let us give the proof, which follows in a couple of lines using the Cauchy-Schwarz inequality: for each $x \in G$ we have
\begin{equation}\label{eq2} \sum_{y \in G} f(y)f(yx) \leq (\sum_{y \in G} f(y)^2 )^{1/2} (\sum_{y \in G} f(yx)^2 )^{1/2} = \sum_{y \in G} f(y)^2,\end{equation} and thus
\begin{equation}\label{eq3} (\sum_{y \in G} f(y) f(yx))^2 \leq (\sum_{y \in G} f(y)f(yx)) \sum_{y \in G} f(y)^2.\end{equation}
Summing over $x \in G$ and using the fact that $\sum_{t \in G} f(t) = 1$, we obtain the result. 

Let us address the equality question, that is to say let us characterise those $f \in \mathscr{F}$ for which $v(f) = 1$. For this to happen, we must have equality in \eqref{eq3} for every $x$. For a given $x$ this means that either $\sum_y f(y) f(yx) = 0$, or else equality occurs in \eqref{eq2}. The first case implies that for all $y$ at least one of $f(y)$ and $f(yx)$ is zero. The second case may be analysed using the well-known criterion for equality in the Cauchy-Schwarz inequality. This implies that there is some $\lambda(x)$ such that $f(y) = \lambda(x)f(yx)$ for all $y$; using the fact that $\sum_{t \in G} f(t) = 1$, it follows that $\lambda(x) = 1$ and therefore $f(y) = f(yx)$ for all $y$. 

Thus $v(f) = 1$ if and only if for all $x \in G$ we have one of the following two mutually exclusive options:
\begin{enumerate}
\item For all $y \in G$, either $f(y)$ or $f(yx)$ is zero;
\item For all $y \in G$, $f(y) = f(yx)$. 
\end{enumerate}
It follows immediately from this that $f$ cannot take more than one non-zero value, and therefore $f(x) = \frac{1}{|A|}1_A(x)$ for some (finite) symmetric set $A \subset G$ containing the identity.  The above two properties then tell us that for all $x \in G$ we have one of the following two mutually exclusive options:
\begin{enumerate}
\item $A$ and $Ax$ are disjoint;
\item $A = Ax$. 
\end{enumerate}
The set of $x$ for which (2) is satisfied is a subgroup of $G$ (the stabiliser of $A$ when $G$ acts on finite subsets of itself by right multiplication). Call this group $H$. If $x \in A$ then, since $\id_G \in A$, the sets $A$ and $Ax$ are not disjoint and so $x \in H$; thus $A \subset H$. On the other hand if $x \in H$ then $A = Ax$ and so in particular, since $\id_G \in A$, we have $x \in A$ and so $H \subset A$. It follows that $A = H$ is a subgroup of $G$. Observations equivalent to these may be found in Hardy-Littlewood \cite{hardy-littlewood}.

Now let us think about the stability and robustness questions. To do this, let us introduce a parameter $K \geq 1$, and let us ask what may be said about those $f \in \mathscr{F}$ for which $v(f) \geq \frac{1}{K}$. This includes both the stability question (where $K \approx 1$) and the robustness question (where $K$ is somewhat larger, for example $K \sim 100$). To spell it out, we are asking for a description of the finitely-supported, symmetric probability measures $f : G \rightarrow [0,\infty)$ for which 
\begin{equation}\label{eq4} \Vert f \ast f \Vert_2^2 \geq \frac{1}{K} \Vert f \Vert_2^2.\end{equation}
To get a feel for this question, let us specialise to the case $f(x) = \frac{1}{|A|} 1_A(x)$, for some finite, symmetric set $A \subset G$ containing the identity. We saw above that only this case is relevant for discussion of the equality question, and in fact the analysis of the stability and robustness questions may be reduced to this case by fairly routine technical arguments \cite{bourgain-gamburd}, \cite[Appendix A]{bggt}. In this case one may check that $\Vert f \Vert_2^2 = |A|^{-1}$ and
\[ \Vert f \ast f \Vert_2^2 = |A|^{-4} \#\{  (a_1,a_2, a_3, a_4) \in A \times A \times A \times A : a_1a_2 = a_3 a_4\}.\]
Thus \eqref{eq4} holds if and only if we have
\begin{equation}\label{eq5} |A|^{-3}\#\{  (a_1,a_2, a_3, a_4) \in A \times A \times A \times A : a_1a_2 = a_3 a_4\} \geq \frac{1}{K} .\end{equation} The quantity on the left here is usually called the \emph{multiplicative energy} $E(A)$ of the set $A$. As can be seen, it records coincidences amongst products of elements of $A$.  Young's inequality implies that $E(A) \leq 1$, and we showed above that equality occurs if and only if $A$ is a subgroup. That $E(A) \leq 1$ can in fact be established easily and directly by noting that if $a_1 a_2 = a_3 a_4$ then $a_4$ is uniquely determined by $a_1,a_2$ and $a_3$. 

When, then, does \eqref{eq5} hold? Here we split the discussion of the stability question ($K \approx 1$) and the robustness question ($K \gg 1$), making just a few remarks about the former. In the stability case it turns out that $A$ must be ``almost'' a subgroup; in fact there is a subgroup $H$ such that the symmetric difference of $A$ and $H$ is very small. Results of this type are certainly very interesting and may be dated to work of Freiman \cite{freiman-groups} and Fournier \cite{fournier} amongst others. Among the diverse applications are the analysis of certain algorithms for sampling at random from finite groups \cite{cooperman, dixon, green-barbados} and the solution of the Dirac-Motzkin conjecture in combinatorial geometry connected with point-line configuations having few ordinary lines \cite{green-tao-pointslines}.

Our main focus here, however, is on the robustness regime $K \gg 1$, where the flavour and the applications are somewhat different. We begin by observing that \eqref{eq5} is implied by a condition which is perhaps easier to understand, that of \emph{small doubling}. We say that a set $A \subset G$ has doubling at most $K$ if
\begin{equation}\label{eq6} |A^2| \leq K|A|,\end{equation}
where $A^2 = \{a_1a_2 : a_1,a_2 \in A\}$. To see that \eqref{eq6} implies \eqref{eq5}, write $r(x)$ for the number of representations of pairs $(a_1, a_2) \in A \times A$ with $a_1 a_2 = x$. Then $r(x) = 0$ for $x \notin A^2$ and so by the Cauchy-Schwarz inequality we have
\[ E(A) = \sum_x r(x)^2 \geq \frac{1}{|A^2|} \big(\sum_x r(x)\big)^2 = \frac{|A|^4}{|A^2|} \geq \frac{1}{K} |A|^3.\]
We have shown that \eqref{eq6} implies \eqref{eq5}, and so if $f(x) = \frac{1}{|A|} 1_A(x)$ for a symmetric set $A$ satisfying \eqref{eq6} then indeed $v(f) \geq \frac{1}{K}$. Thus an analysis of the robustness question for Young's inequality necessarily involves studying sets $A$ satisfying \eqref{eq6}. It is not at all obvious that such a study is sufficient for that task, because we have not shown that \eqref{eq5} implies \eqref{eq6}. In fact, it does not, as be easily seen by taking $A$ to be $H \cup X$, where $H$ is a subgroup of $G$ and $X$ is an arbitrary symmetric set of the same size, disjoint from $H$. Then \eqref{eq5} holds with $K = 8$, since we may take all quadruples $(a_1,a_2,a_3,a_4)$ with $a_1,\dots, a_4 \in H$ and $a_1a_2 = a_3 a_4$. However, there is absolutely no reason to suppose that \eqref{eq6} holds, and indeed $A^2$ contains $X^2$ which could have size as large as $c |X|^2$. We leave it to the reader to provide an explicit example in a suitable group $G$. Remarkably, however, the large multiplicative energy condition \eqref{eq5} does imply a weak version of \eqref{eq6}: specifically, \eqref{eq6} is true after passing from $A$ to a large subset $A'$ and replacing \eqref{eq6} by a somewhat weaker condition $|A^{\prime 2}| \leq K'|A'|$ with $K' \sim K^{10}$, say. This result is known as the Balog-Szemer\'edi-Gowers theorem, because in the case $G$ abelian it was established by Gowers \cite{gowers-4aps} in the course of his seminal work on Szemer\'edi's theorem, an earlier result of a qualitatively similar form but with the bound on $K'$ being vastly weaker having previously been established by Balog and Szemer\'edi \cite{bs} by different means. It was shown by Tao \cite{tao-noncommutative} that the assumption that $G$ is abelian could be dropped. 

\subsection{Approximate groups} 
We have discussed the relationship between the robustness question for Young's inequality and the study of finite sets $A$ satisfying the small doubling condition $|A^2| \leq K|A|$. Since subgroups of $G$ provide equality in Young's inequality, this provides some justification for thinking of such $A$ as ``approximate groups''. Moreover, the small doubling condition visibly suggests that $A$ is somehow almost closed under multiplication, surely a property we would expect from any sensible notion of an approximate group. As it turns out, it has been found convenient to introduce a slightly different but closely related notion.

\begin{definition}[Approximate group]
Let $A$ be a subset of a group $G$. Then we say that $A$ is a $K$-approximate group if $A$ is symmetric, contains the identity, and if $A^2 \subset XA$ for some set $X$ of size at most $K$. 
\end{definition}
This definition was introduced by Tao \cite{tao-noncommutative} and has certain advantages such as behaving well under homomorphisms, making sense for infinite sets $A$ as well as finite ones, and immediately implying further conditions on $A$ such the tripling bound $|A^3| \leq K^2 |A|$.

Note that if $A$ is a $K$-approximate group then $A$ automatically satisfies the small doubling condition \eqref{eq6}, and hence the large multiplicative energy condition \eqref{eq5}. The reverse direction is less clear, and the situation is much the same as before: a set satisfying \eqref{eq6} need not be a $K$-approximate group, but there is a closely related set $A'$ which is a $K'$-approximate group for some $K' \sim K^{10}$. This deduction is essentially due to Ruzsa, who laid the foundation for the whole theory in a series of works.  For a precise statement and further references, \S 4 of \cite{bgt-survey} may be consulted. 
\subsection{Examples}
We now give some examples of approximate groups. The first example is fairly trivial.

\begin{example}\label{ex55}
If $A \subset G$ is a subgroup then of course $A$ is symmetric, $\id_G \in A$ and $A^2 = A$. Thus $A$ is a $1$-approximate group. 
\end{example}

Thus far, we have not pointed out that there are in fact nontrivial examples of approximate groups. The simplest is a geometric progression.

\begin{example}\label{ex66}
If $P$ is the geometric progression \[ P = P(u; N) := \{u^{n} : 0 \leq n< N\}\] for some element $u \in G$ and if $A = P \cup P^{-1}$ then $A$ is a $2$-approximate group. Indeed $A = \{u^{n} : -N+1 \leq n \leq N-1\}$, $A^2 = \{u^{n} :  -2N + 2 \leq  n \leq 2N - 2\}$ and so $A^2 \subset X A$ where $X = \{u^{N - 1}, u^{-N + 1}\}$. 
\end{example}

Less obviously, there are multidimensional generalisations of the preceding example. 

\begin{example}\label{ex7}
If $P$ is the multidimensional geometric progression 
$$P = P(u_1,\dots, u_d; N_1, \dots, N_d) := \{u_1^{n_1} u_2^{n_2} \dots u_d^{n_d} : 0 \leq n_i < N_i\}$$ for some \emph{commuting} elements $u_1,\dots, u_d \in G$ and integers $N_1,\ldots,N_d>0$ and if $A = P \cup P^{-1}$ then $A$ is a $2^d$-approximate group. We leave the confirmation of this to the reader.
\end{example}

The commuting assumption was very important in the previous example (otherwise we cannot simplify a product $u_1^{n_1} \dots u_d^{n_d} u_1^{n'_1} \dots u_d^{n'_d}$ to $u_1^{n_1 + n'_1} \dots u_d^{n_d + n'_d}$). However, it can be replaced by the weaker condition of nilpotence, as the following example shows. 

\begin{example}\label{ex8} Let $N_1,N_2,N_{1,2}$ be positive integers with $N_{1,2} \geq N_1 N_2$, let $G = \heis{\R}{\R}{\R}$ be the Heisenberg group and let $A \subset G$ be the following set of matrices. Let
\[ u_1 :=  \heis{1}{0}{0} , \quad u_2 := \heis{0}{1}{0} ,\] and take $A = P \cup P^{-1}$ where $P = P(u_1, u_2, [u_1,u_2]; N_1, N_2, N_{1,2})$ is the set
\[ \{ u_1^{n_1} u_2^{n_2} [u_1, u_2]^{n_{1,2}} : 0 \leq n_1 < N_1, 0 \leq n_2 < N_2, 0 \leq n_{1,2} < N_{1,2}\}.
\]
Here, $[u_1,u_2]$ is the commutator given by
\[ [u_1, u_2] := u_1 u_2 u_1^{-1} u_2^{-1} = \heis{0}{0}{1}.\] 
It may be straightforwardly checked that 
\begin{equation}\label{prod} u_1^{n_1} u_2^{n_2} [u_1, u_2]^{n_{1,2}} \cdot u_1^{n'_1} u_2^{n'_2} [u_1, u_2]^{n'_{1,2}} = u_1^{n_1 + n'_1} u_2^{n_2 + n'_2} [u_1, u_2]^{n_{1,2} + n'_{1,2} - n'_1 n_2}.\end{equation}
Hence
\[ P^{-1} \subset \{ u_1^{n_1} u_2^{n_2} [u_1, u_2]^{n_{1,2}} : -N_1 < n_1 \leq 0, -N_2 < n_2 \leq 0, - 2N_{1,2} < n_{1,2}\leq 0\} \]
and
\[ A^2 \subset \{ u_1^{n_1} u_2^{n_2} [u_1, u_2]^{n_{1,2}} : |n_1| < 2N_1, |n_2| < 2N_2, |n_{1,2}| < 5N_{1,2}\}. \]
Now for any $n'_1, n'_2, n'_{1,2}$ in \eqref{prod} we may choose (unique) integers $k_1, k_2, k_{1,2}$ such that 
\[ u_1^{k_1 N_1} u_2^{k_2 N_2} [u_1, u_2]^{k_{1,2} N_{1,2}} \cdot u_1^{n'_1} u_2^{n'_2} [u_1, u_2]^{n'_{1,2}}  \in P.\] Indeed we have $k_1 = -\lfloor n'_1/N_1\rfloor$, $k_2 = -\lfloor n'_2/N_2\rfloor$ and $k_{1,2} = -\lfloor (n'_{1,2} - n'_1 k_2 N_2)/N'_{1,2} \rfloor$. Thus if $u_1^{n'_1} u_2^{n'_2} [u_1, u_2]^{n'_{1,2}} \in A^2$ then $|k_1| \leq 1$, $|k_2| \leq 1$ and $|k_{1,2}| \leq 6$. Hence
\[ A^2 \subset XP \subset XA,\] where
\[ X = \{ u_1^{k_1 N_1} u_2^{k_2 N_2} [u_1, u_2]^{k_{1,2} N_{1,2}} : |k_1| \leq 1, |k_2| \leq 1, |k_{1,2}| \leq 6\} \] is a set of size 117. That is, $A$ is a $117$-approximate group. (A smaller constant could be obtained with a more careful analysis.)
\end{example}

Example \ref{ex8} is an example of a \emph{nilprogression}. The key feature of the Heisenberg group $G$ relevant to this example is the fact that it is nilpotent of class $2$, which means that commutators of order 3 or higher are all equal to the identity, or equivalently that $[u_1, u_2]$ commutes with everything else. Similar examples can be constructed in more general nilpotent groups of arbitrary class $s$, though the constant $K$ (117 in Example \ref{ex8}) will generally grow with $s$.  We will not give the details here, and refer the reader instead to \cite[Definition 2.1]{bgt-survey}. The nilprogression in Example \ref{ex8} is said to have rank $2$ and class $2$ (the rank being the number of generators $u_i$ and the class being the nilpotency class of the group generated by $u_1,u_2$). 

Different instances of the above constructions may be combined to create new examples. For example, it is easy to see that the direct product of a $K_1$-approximate group and a $K_2$-approxiate group is a $(K_1K_2)$-approximate group. There are also combinations of the above examples which are not direct products, for example the following example of Helfgott \cite{helfgott-sl3}.

\begin{example}\label{ex9}
Let $p$ be a large prime, let $r, s, t \in \F_p$ be fixed generators of $\F_p^*$, let $N_1,N_2,N_3$ be positive integers, and define $A$ to be a set of $3 \times 3$ matrices over $\F_p$ as follows. Set $A = H P \cup (H P)^{-1}$, where 
\[ H := \left\{ \left( \begin{smallmatrix} 1 & x& z \\ 0 & 1 & y \\ 0 & 0 & 1 \end{smallmatrix} \right) : x, y, z \in \F_p\right\}\]
and
\[ P = P(u_1, u_2, u_3; N_1, N_2, N_3) :=\{ u_1^{n_1} u_2^{n_2} u_3^{n_3} : 0 \leq n_i < N_i\}\] 
with 
\[ u_1 := \left( \begin{smallmatrix} r & 0 & 0 \\ 0 & 1 & 0 \\ 0 & 0 & 1 \end{smallmatrix} \right), u_2 := \left( \begin{smallmatrix} 1 & 0 & 0 \\ 0 & s & 0 \\ 0 & 0 & 1 \end{smallmatrix} \right),  u_3 := \left( \begin{smallmatrix} 1 & 0 & 0 \\ 0 & 1 & 0 \\ 0 & 0 & t \end{smallmatrix} \right)\] for some $r, s, t \in \F_p^*$.
as in Example \ref{ex7} above. It is quite easy to check that 
\[ A = \left\{ \left( \begin{smallmatrix} r^{n_1} & x& z \\ 0 & s^{n_2} & y \\ 0 & 0 & t^{n_3} \end{smallmatrix} \right) : x,y,z \in \F_p, -N_i < n_i < N_i\right\}\]
and hence
\[ A^2 \subset \left\{ \left( \begin{smallmatrix} r^{n'_1} & x& z \\ 0 & s^{n'_2} & y \\ 0 & 0 & t^{n'_3} \end{smallmatrix} \right) : x,y,z \in \F_p, -2N_i  < n'_i < 2N_i\right\},\] from which it follows that $A$ is an $8$-approximate group.\end{example}

In Example \ref{ex9}, $H P$ was an example of a \emph{coset nilprogression}, in this case of rank $3$ and step $1$. The general form of a coset progression is $H P$ where $P$ is a nilprogression (a notion we did not define in full generality) and $H$ is a subgroup normal in the group $\langle P \rangle$ generated by $P$. In fact, all five of our examples were of the form $A = (H P) \cup (H  P)^{-1}$ for some coset progression $H  P$ (in Examples \ref{ex66}, \ref{ex7} and \ref{ex8} the subgroup $H$ was trivial, whilst in Example \ref{ex55} the nilprogression $P$ was trivial). Conversely, every $A$ of this form  is a $K$-approximate group, where $K$ is bounded as a function of the rank $r$ and the class $s$ of $P$.  Once again we refer the reader to \cite{bgt-survey} for more information.

\subsection{Theorems about approximate groups} 

Given the discussion of the last section, it is natural to ask whether every $K$-approximate group is of the form $(H P) \cup (H  P)^{-1}$ for some coset nilprogression $H P$ (of rank and step bounded in terms of $K$). The answer to this is, strictly speaking, negative, as the following example of a set $A \subset \Z$ shows. Here, we use additive notation for the group operation on $\Z$ and so our interest is in $2A = A + A$ rather than $A^2$. Define $A$ to be $\{0\} \cup \bigcup_{j = 1}^{N} \{2j - \eps_j, -2j + \eps_j\}$, where the $\eps_j$ are independent $\{0,1\}$-valued random variables. Then $2A  \subset [-4N, 4N]$. However, $\{-1, 0, 1\} + A  \supset [-2N, 2N]$, and so $\{ - 2N , 2N\} + \{-1,0,1\} + A \supset [-4N, 4N] \supset 2A $.  It follows that $A$ is a $6$-approximate group. However, for a typical choice of the $\eps_j$, $A$ does not have nearly so much structure as a progression (though it is \emph{syndetic}, that is to say has bounded gaps, which is what makes this construction work). 

However we do have the following recent result of Breuillard, Tao and the author \cite{bgt}. 

\begin{theorem}\label{bgt-main}
Let $A$ be a $K$-approximate subgroup of a group $G$. Then there is a coset nilprogression $B = H P$ of rank and class bounded as functions of $K$, where such that $|B| \leq K'|A|$ and there is a set $X \subset G$ with $|X| \leq K'$ such that $A \subset (X B) \cap (B X)$. Here, $K'$ may be bounded as a function of $K$ only.
\end{theorem}

We say that $A$ is $K'$-controlled by $B$. In the example preceding the theorem, we may take $B = \{0,\dots, N-1\}$. The reader is encouraged not to dwell too lengthily on the notion of ``control'' and read the above theorem as follows: every approximate group is roughly a coset nilprogression.

\begin{itemize}
\item For many specific types of group $G$, statements equivalent to Theorem \ref{bgt-main} had previously been established, often with good quantitative control over the parameter $K'$ as well as the rank and class. When $G = \Z$, this is essentially the celebrated Freiman-Ruzsa theorem \cite{freiman-book, ruzsa-freiman}. The general abelian case was handled by Ruzsa and the author \cite{green-ruzsa}, building on earlier work of Ruzsa \cite{ruz-fin-field}. Various matrix groups $G$ were handled in work of (in chronological order) Elekes-Kir\'aly \cite{ek}, Chang \cite{chang} and Helfgott \cite{helfgott-sl2,helfgott-sl3}, the latter handling $\mbox{SL}_2(k)$ and $\mbox{SL}_3(k)$ with $k = \F_p$ or $k = \C$, amongst others.
\item Hrushovski \cite{hrush}, in a very important 2009 breakthrough, dealt with $G = \mbox{GL}_n(\C)$ (though with some dependence on $n$). His argument was model-theoretic and a key ingredient of it was his ``Lie model theorem'', also a key ingredient in the proof of Theorem \ref{bgt-main}.
\item The proof of Theorem \ref{bgt-main} additionally requires arguments related to the solution of Hilbert's Fifth Problem (every locally compact group is locally an inverse limit of Lie groups), specifically lemmas due to Gleason from the 1950s. It also makes use of a lemma in additive combinatorics of a type developed by Sanders \cite{sanders-lemma} and Croot-Sisask \cite{croot-sisask}. 
\item Theorem \ref{bgt-main} is in fact valid when $G$ is a ``local group'' rather than a \emph{bona fide} group. Moreover, it was necessary in \cite{bgt} to work in this larger category, although Hrushovski and van den Dries have since managed to arrange the argument so that, at the expense of proving a slightly weaker result,  one need only work in genuine groups.
\item Theorem \ref{bgt-main} rather easily implies Gromov's famous theorem \cite{gromov} on groups of polynomial growth. However, it does not really provide a new proof of Gromov's theorem as all the deep ingredients Gromov developed (the notion of an asymptotic cone, and the application of the solution to Hilbert's Fifth Problem) are also required here in some form. 
\end{itemize}

Whilst Theorem \ref{bgt-main} is definitive from the qualitative point of view, for many applications more quantitative statements are required. Unfortunately, a crucial use of ultrafilters in the proof of Theorem \ref{bgt-main} means that no quantitative dependence of $K'$ on $K$ is currently known\footnote{In principle one could be obtained by quantifier elimination but this would be a huge amount of effort and the bound would be desperately weak.}. In the next section we will discuss perhaps the most substantial application of the theory of approximate groups so far, to the study of rapidly mixing random walks on groups (expanders), which find further application in the ``affine sieve''. For these applications much more quantitative statements are required, but in more restricted settings.

A celebrated result of the type we have in mind is the theorem of Helfgott \cite{helfgott-sl2}.

\begin{theorem}[Helfgott]\label{helfgott-thm}
Let $K \geq 2$. Suppose that $A \subset G$ is a $K$-approximate group, where $G = \PSL_2(\F_p)$, and that $A$ generates $G$. Then either $|A| \leq K^{C_2}$ or $|A| \geq K^{-C_2} |G|$, where $C_2$ is an absolute constant. 
\end{theorem}

This theorem was generalised to $\PSL_3(\F_p)$ by Helfgott in a subsequent paper \cite{helfgott-sl3}, and then to $\PSL_n(\F_p)$ (with $C_2$ replaced by an exponent $C_n$ depending on $n$) and other finite simple groups of Lie type in independent works of Pyber-Szab\'o \cite{pyber-szabo} and Breuillard, Tao and the author \cite{bgt-linear}, the former paper containing a slightly more general result than the latter.

It is worth remarking that Helfgott's arguments made substantial use of the theory of \emph{approximate fields} or ``sum-product theory'', in particular a result of Bourgain, Katz and Tao \cite{bkt}. This is an important topic in arithmetic combinatorics and it has links to the theory of approximate groups as well as other substantial applications, perhaps most notably estimates for the additive Fourier transform of multiplicative subgroups of $\F_p^{\times}$ due to Bourgain, Glibichuk and Konyagin \cite{bgk}. The subsequent works \cite{bgt-linear,pyber-szabo} do not make explicit use of this theory, and in fact it was noted in \cite{bgt-linear} that, conversely, results about approximate subgroups of $\SL_2(k)$ imply results about approximate subfields of $k$. Sadly we do not have the space to discuss these aspects any further here.

It is also of interest to note that \cite{bgt-linear} made use of an analogue for approximate groups of an argument of Larson and Pink \cite{larson-pink}, which gives a self-contained and relatively concise proof of certain statements which follow from the Classification of Finite Simple Groups (CFSG).

\subsection{Applications}

Several applications of Theorem \ref{bgt-main} are given in the paper \cite{bgt}. In addition to certain refinements of Gromov's theorem they include a result about the virtual nilpotence of the fundamental group of almost negatively-curved Riemannian manifolds, and a generalisation of a lemma of Margulis stating that the ``almost stabiliser'' of a point $x$ in a finite-dimensional metric space $X$ under the action of a discrete group of isometries is virtually nilpotent.  Here, however, we wish to discuss an appealing application, due to Bourgain and Gamburd \cite{bourgain-gamburd} of Helfgott's result, Theorem \ref{helfgott-thm}.

The result concerns a property of a generating set $S$ of a finite group $G$ known as \emph{expansion}. This property has several equivalent characterisations, details of which may be found in the survey \cite{hoory-etc}. For our purposes here, however, it is convenient to define expansion in terms of the rapid mixing of the random walk on generators $S \cup S^{-1}$. For the sake of illustration suppose that $|S| = 2$, write $S= \{a,b\}$, and imagine $G$ being quite large. Then we perform a random walk of $m$ steps, the end result of which is a product $x_m = g_1 \dots g_m$ where each $g_i$ is selected independently at random from the set $\{a, b, a^{-1}, b^{-1}\}$.  Note that if $G = \Z^2$ and $S = \{(1,0), (0,1)\}$ (and if additive notation is used) then this is precisely the classical random walk on the plane $\R^2$.

Now $x_1$ takes values in a set of size $4$, tiny in comparison to $|G|$, and by a trivial induction $x_j$ takes values in a set of size at most $4^j$ (in fact by an almost-as-trivial induction one may reduce this to $4 \cdot 3^{j-1}$). Thus if $j = c \log |G|$ for some small value of $c$ then $x_j$ takes values in a set of size at most $|G|^{c'}$, and in particular is nowhere near to equidistributed on $G$. However in certain situations it turns out to be the case that $x_j$ \emph{is} highly equidistributed not much later than this time, say for $j \geq C \log |G|$, for some $C$. By ``highly-equidistributed'' let us (slightly arbitrarily) say that we mean $\mathbb{P}(x_j = g) = \frac{1}{|G|} + O(\frac{1}{|G|^{10}})$ for all $g \in G$. The situation just described is one possible definition of what it means for $S$ to be an expander (the precise definition must include the parameter $C$). 

\begin{theorem}[Bourgain-Gamburd]\label{bourg-gamb}
Let $G = \PSL_2(\F_p)$ and suppose that $S = \{ \left(\begin{smallmatrix} 1 & 3 \\ 0 & 1 \end{smallmatrix}\right), \left(\begin{smallmatrix} 1 & 0 \\ 3 & 1 \end{smallmatrix}\right)\}$. Then the random walk on generating set $S \cup S^{-1}$ becomes highly equidistributed in time at most $C \log |G|$ for some absolute constant $C$,  independent of $p$. 
\end{theorem}

The only reason we have put ``3'' in the matrices here is that the result was already known with ``1'' and ``2'' by different methods. Bourgain and Gamburd actually proved a far more general result, but we fix on this special case for the sake of illustration. To describe the proof, we note that the result can be reformulated in terms of convolution powers of $\mu = \frac{1}{4}(\delta_a + \delta_b + \delta_{a^{-1}} + \delta_{b^{-1}})$, that is to say the function $\mu :G \rightarrow [0,\infty)$ taking values $\frac{1}{4}$ at $x = a, b , a^{-1}, b^{-1}$ and zero elsewhere. The probability that $x_j = g$ is then precisely $\mu^{(j)}(g)$, where $\mu^{(j)} = \mu \ast \dots \ast \mu$ and the convolution is repeated $j$ times, and where we are defining $\nu_1 \ast \nu_2 (x) = \sum_y \nu_1(y) \nu_2(x y^{-1})$. Note that $\sum_x \mu^{(j)}(x) = 1$. 

We are interested in how quickly $\mu^{(j)}$ tends towards the constant function $\frac{1}{|G|}$. To study this, we follow the progress of $\mu^{(j)}$ in three stages:

\begin{itemize}
\item The early stage in which $j \leq \frac{1}{10}\log |G|$;
\item The middle stage in which $\frac{1}{10} \log |G| \leq j \leq \frac{C}{10}\log |G|$;
\item The end stage in which $\frac{C}{10}\log |G| \leq j \leq C \log|G|$.
\end{itemize}
The early stage is relatively easy to analyse. This is because the elements $a, b$ behave as though they were generators of a free group, or in other words the random walk does not intersect itself nontrivially, so long as words of length at most $\frac{1}{10} \log |G|$ are being considered. As a consequence, the size $|\mbox{Supp}(\mu^{(j)})|$ of the support of $\mu^{(j)}$ at the end of the early stage is somewhat large, of size at least about $|G|^{0.01}$, say. In fact the support is not quite the most sensible thing to look at, because $\mu^{(j)}$ may take on several different values. A more nuanced quantity is $\Vert \mu^{(j)}\Vert_2^{-2}$, which we will call the \emph{weighted support}. This would equal $|\mbox{Supp}(\mu^{(j)})|$ if $\mu^{(j)}$ did happen to be constant on its support, as is easily checked. 

The theory of approximate groups is applied to analyse the middle stage. If $j_1$ is the end of the early stage, we look at iterates $\mu^{(j_1)}, \mu^{(2j_1)}, \dots, \mu^{(j_2)}$ where $j_2= 2^{\ell} j_1$ for some $\ell$. If $\ell$ is somewhat large, the weighted support of $\mu^{(2^t j_1)}$ cannot always increase substantially as we change $t$ to $t+1$, and so there must be some $t$ for which the weighted supports of $\mu^{(2^t j_1)}$ and of $\mu^{(2^{t+1} j_1)}$ are roughly the same. Since $\mu^{(2^{t+1} j_1)} = \mu^{(2^t j_1)} \ast \mu^{(2^t j_1)}$, the only way that this can happen is if $f = \mu^{(2^t j_1)}$ satisfies \eqref{eq4} for some fairly small value of $K$, that is to say $\Vert f \ast f \Vert_2^2 \geq \frac{1}{K} \Vert f \Vert_2^2$ or $v(f) \geq \frac{1}{K}$. This is precisely the robustness question for Young's inequality that we have been studying. As we discussed, this situation implies, very roughly speaking\footnote{The $\sim$ notation here hides quite a few technicalities.}, that $f \sim \frac{1}{|A|}1_A(x)$ where $A$ is a $K$-approximate group. Here of course we are concerned with the particular case $G = \PSL_2(\F_p)$, so by Helfgott's Theorem \ref{helfgott-thm} there are three possibilities: (i) $A$ is tiny, (ii) $A$ is almost all of $G$ and (iii) $A$ does not generate $G$. Case (i) cannot occur, because at the end of the early stage the weighted support of $\mu^{(j_1)}$ was quite large. It turns out that (iii) also cannot occur, because of the particular structure of proper subgroups of $\PSL_2(\F_p)$: they are all soluble and so satisfy the law $[[x_1, x_2], [x_3, x_4]] = \id_G$, quite at odds with the free behaviour exhibited during the early stage. We are left, then with possibility (ii), which implies that the weighted support of $\mu^{(2^t j_1)}$ is almost $|G|$. By further applications of Young's inequality the same is true of $\mu^{(2^{\ell} j_1)} = \mu^{(j_2)}$. That is to say, at the end of the middle stage $\mu^{(j_2)}$ fills out a large portion of $G$ in a fairly uniform way.    

The analysis of the end stage involves still different ideas -- an application of representation theory having its origin in a paper of Sarnak and Xue \cite{sarnak-xue}. The crucial input is the fact that all nontrivial representations of $G = \PSL_2(\F_p)$ have dimension at least $\frac{1}{2}(p-1)$ and in particular at least $|G|^c$ for some constant $c$. (In the language of Gowers \cite{gowers-quasirandom},  $G$ is an example of a ``quasirandom'' group.) Further details may, of course, be found in the original paper \cite{bourgain-gamburd}.

The ``Bourgain-Gamburd expansion machine'' just described and modifications of it have found many further applications. One is the following variant of Theorem \ref{bourg-gamb} due to Breuillard, Guralnick, Tao and the author \cite{bggt}.

\begin{theorem}\label{thm2.5}
Let $G$ be any finite simple group of Lie type and suppose that $S = \{ a,b\}$ where $a,b$ are chosen uniformly at random from $G$. Then, with probability at least $1 - O(|G|^{-c})$, the random walk on generating set $S \cup S^{-1}$ becomes highly equidistributed in time at most $C \log |G|$. Here $c, C > 0$ depend only on the rank of $G$. 
\end{theorem}

For example, this theorem holds with $G = \PSL_n(\F_q)$ and with $C$ depending only on $n$ and not on $q$. The proof of this theorem relies on the Bourgain-Gamburd expansion machine but with the work of Pyber--Szab\'o \cite{pyber-szabo} and Breuillard, Tao and the author \cite{bgt-linear} in place of Helfgott's theorem. It also requires several other ingredients, including two different \emph{ad hoc} analyses in two particular families of groups (the symplectic groups $\mbox{Sp}_4(k)$ in characteristic $3$ and the triality Groups ${}^3D_4(q)$). A different particular case, that in which $G$ is a Suzuki group $\mbox{Sz}(q)$, had been handled in an earlier paper \cite{bgt-suzuki} of the authors. This was of a certain amount of interest because it completed the proof of the following theorem of Lubotzky, Kassabov and Nikolov \cite{lkn}.

\begin{theorem}
There are absolute constants $k, C$ with the following property. For any nonabelian finite simple group $G$, there is a set $S \subset G$ of size at most $k$ such that the random walk on generating set $S \cup S^{-1}$ becomes highly equidistributed in time at most $C \log |G|$.
\end{theorem}

The proof of this theorem depends on CFSG and the most impressive ingredient is, in my view, Kassabov's proof \cite{kassabov} in the case $G = A_n$. It appears to be unknown whether or not Theorem \ref{thm2.5} holds uniformly for all finite simple groups $G$ with $C$ an absolute constant, even in (especially in?) the case $G = A_n$. 

Perhaps of greater interest for applications than results such as Theorem \ref{thm2.5}, however, are generalisations of the original Bourgain-Gamburd theorem, where the groups under consideration range over a family such as $G = \PSL_n(\F_p)$, $p$ prime and the set $S$ is obtained by reduction of a \emph{fixed} set of integer matrices, rather than by random selection for each $p$. In the Bourgain-Gamburd theorem as stated above, $S =  \{ \left(\begin{smallmatrix} 1 & 3 \\ 0 & 1 \end{smallmatrix}\right), \left(\begin{smallmatrix} 1 & 0 \\ 3 & 1 \end{smallmatrix}\right)\}$. The crucial property of this set of generators for rapid mixing of the random walk is that, considered as a subset of $\mbox{SL}_2(\Z)$, the subgroup they generate is \emph{Zariski dense} (not contained in any proper algebraic subvariety). That this condition is sufficient was established by Bourgain-Gamburd for the family $\PSL_2(\F_p)$, $p$ prime.  Varj\'u \cite{varju} obtained the same result for $\PSL_n(\F_p)$, and moreover for $\PSL_n(\Z/q\Z)$ where $q$ is squarefree but may well be composite. (Such results had already been established in the case $n = 2$ by Bourgain, Gamburd and Sarnak \cite{bgs-sl2-fq} by a more complicated method based in part on Helfgott's arguments, necessitating in particular a foray into the tricky territory of \emph{approximate subrings} of $\Z/q\Z$.) This last result is a crucial ingredient in the so-called \emph{affine sieve} of Bourgain, Gamburd and Sarnak which finds almost primes in the matrix entries of orbits in matrix groups. Any serious discussion of this would take us too far afield, so we refer the reader to \cite{sarnak-salehi} for the state of the art and to the very nice exposition \cite{sarnak-notes} for a (somewhat outdated) introduction. See also \cite{green-ceb}, again rather outdated.

\subsection{Open questions}

There are many open questions concerning the quantitative aspects of the theory described above. For example, no version of Theorem \ref{bgt-main} in which the parameter $K'$ is given quantitatively in terms of $K$ is known, and nor does it seem prudent at this stage to speculate on what might be true in this regard. One tempting line of enquiry would be to look at Kleiner's alternative proof \cite{kleiner} of Gromov's theorem in the context of approximate groups, but this has not so far been successful. 

Even in the case $G = \Z$ there are unsolved problems connected with approximate groups. As previously noted, Theorem \ref{bgt-main} in this case is due to Freiman \cite{freiman-book} and Ruzsa \cite{ruzsa-freiman}. In $\Z$, there are no interesting finite subgroups and, of course, all nilprogressions are automatically abelian progressions as in Example \ref{ex7}. Writing the group operation on $\Z$ using addition as usual, the Freiman-Ruzsa theorem may be stated as follows.

\begin{theorem}[Freiman-Ruzsa]\label{freiman-ruzsa-theorem}
Suppose that $A \subset \Z$ is a $K$-approximate group, that is to say $2A \subset X + A$ for some set $X \subset \Z$ with $|X| \leq K$. Then there is a proper\footnote{This means that all the sums $n_1 u_1 + \dots + n_d u_d$ under consideration are distinct.} progression $P = P(u_1,\dots, u_d; N_1, \dots, N_d) := \{ n_1 u_1 + \dots + n_d u_d : 0 \leq n < N\}$ which $K'$-controls $A$. Here, $d$ and $K'$ are bounded as functions of $K$ only. 
\end{theorem}

The definition of ``control'' here is the same as in Theorem \ref{bgt-main}. 

The optimal bounds on $d$ and $K'$ are not known. Following a sequence of developments by Chang \cite{chang-freiman} and Schoen \cite{schoen-freiman}, the state of the art is contained in a breakthrough paper of Sanders \cite{sanders-bogolyubov}. Sanders shows that we may take $d \sim (\log K)^C$ and $K' \sim e^{(\log K)^C}$ for some reasonable value of $C$ (such as $C = 4$). A key open question, known as the \emph{Polynomial Freiman-Ruzsa conjecture}, asks whether one could in fact take $d \sim \log K$ and $K' \sim K^C$. The bound $d \sim \log K$ is significant as if $P$ is a progression of this dimension then $P \cup -P$ is itself a $K^{C'}$-approximate group. If one is prepared to sacrifice $K'$ then bounds of this strength are known due to work of Freiman-Bilu \cite{bilu} and Tao and the author \cite{gt-freiman-bilu}. For much greater depth on the quantitative issues surrounding Theorem \ref{freiman-ruzsa-theorem}, the recent survey of Sanders \cite{sanders-survey} may be consulted.

A solution to the Polynomial Freiman-Ruzsa conjecture ought to have serious applications in additive number theory -- perhaps, for example, to questions about bases such as Waring's problem. However, no definite deductions of this type have so far been made. 

Another abelian setting has attracted a lot of interest, and that is the case $G = \F_2^{\Z}$. In this group, where we have $2 \cdot x = 0$ for every $x$, there are no interesting nilprogressions and one is left only with subgroups. Theorem \ref{bgt-main} in this case is due to Ruzsa \cite{ruz-fin-field}, and it may be stated as follows.

\begin{theorem}[Ruzsa]\label{ruz-fin-thm}
Suppose that $A \subset \F_2^\Z$ is a $K$-approximate group, that is to say $2A  \subset X + A$ for some set $X \subset \F_2^\Z$ with $|X| \leq K$. Then there is a subgroup $H \subset \F_2^{\Z}$ which $K'$-controls $A$. Here $K'$ is bounded as a function of $K$ only. 
\end{theorem}

The question of whether $K'$ may be taken to be polynomial in $K$ is also known as the Polynomial Freiman-Ruzsa conjecture, and it has attracted much attention. Ruzsa \cite{ruz-fin-field} attributes it to Katalin Marton. Once again the best results are due to Sanders \cite{sanders-bogolyubov}, who shows that we may take $K' \sim e^{(\log K)^C}$. Ruzsa (unpublished, but see \cite{green-ruzsa-notes}) offers several equivalent formulations, of which the following is perhaps particularly appealing.

\begin{conjecture}\label{conj2.9}
Let $V$ be a finite-dimensional vector space in characteristic $2$. Suppose that $f : V \rightarrow V$ satisfies the ``approximate homomorphism'' condition
\[ \{ f(x + y) - f(x) - f(y) : x, y \in V \} \subset S.\]
Then there is a linear map $\tilde f : V \rightarrow V$ and a set $\tilde S$ with $|\tilde S| \ll |S|^C$ such that 
\[ \{ f(x) - \tilde f (x) : x \in V\} \subset \tilde S.\]
\end{conjecture}
There is an extremely extensive literature on the closely-related notion of a \emph{quasimorphism} in contexts arising in geometric group theory; see \cite{quasimorphisms} for a brief introduction. At present there seems to be little connection between that context, where the concern is usually with quasimorphisms on infinite groups, and ours.

\section{Approximate polynomials}

\subsection{Gowers norms and polynomial phases} We turn now to the discussion of a different inequality. If $f : \Z \rightarrow \C$ is a function and $h \in \Z$ then we define the multiplicative derivative $\Delta_h f$ by $\Delta_h f(x) = f(x) \overline{f(x + h)}$. Let $k \geq 2$ be a fixed integer, and suppose that $N$ is large in terms of $k$. Write $[N] = \{1,\dots, N\}$. Then we define the \emph{Gowers $U^k[N]$-norm} of $f$ by 
\[ \Vert f \Vert_{U^k[N]} =  \big(\E_{x, h_1, \dots, h_k} \Delta_{h_1} \dots \Delta_{h_k} f(x)\big)^{1/2^k}.\]
Here, the average $\E$ is over all $x, h_1, \dots, h_k$ for which $x + \omega_1 h_1 + \dots + \omega_k h_k \in [N]$ for all $\omega_i \in \{0,1\}$; this means that the Gowers norm depends only on the values taken by $f$ on $[N]$. In taking $2^k$th roots we make use of the not completely obvious fact that $\E_{x, h_1, \dots, h_k} \Delta_{h_1} \dots \Delta_{h_k} f(x)$ is real and non-negative. This is not too hard to prove by induction: see for example \cite{tao-vu}. The basic theory of Gowers norms was originally developed in \cite{gowers-longaps}.

The Gowers norms satisfy the following rather trivial inequality: if $\mathscr{F}$ is the set of all functions $f : [N] \rightarrow \C$ with $\Vert f \Vert_{\infty} \leq 1$,  $v(f) = \Vert f \Vert_{U_k[N]}$, then
\begin{equation}\label{gowers-triv} v(f) \leq 1.\end{equation}
(The inequality is indeed trivial -- bound every instance of $f( \cdot)$ in the definition of the Gowers $U^k$-norm by $1$).

When does equality occur, that is to say for which $f$ do we have $v(f) = 1$? For this to happen, we must have\footnote{Here and in what follows we ignore the restriction that $x + \omega_1 h_1 + \dots + \omega_k h_k \in [N]$; this has little bearing on the argument.}
\begin{equation}\label{mult-der-cond} \Delta_{h_1} \dots \Delta_{h_k} f(x) = 1 \quad \mbox{for all $x, h_1,\dots, h_k$}.\end{equation} This implies that $|f(x)| = 1$ for all $x$, and so we may write $f(x) = e^{2\pi i \phi(x)}$ for some phase function $\phi : \Z \rightarrow \R/\Z$. The condition \eqref{mult-der-cond} then becomes 
\begin{equation}\label{add-der-cond} \partial_{h_1} \dots \partial_{h_k} \phi(x) = 0 \quad \mbox{for all $x, h_1,\dots, h_k$} ,\end{equation} where $\partial_h$ is the additive difference operator defined by $\partial_h \psi(x) = \psi(x) - \psi(x + h)$. 

The condition \eqref{add-der-cond} is satisfied if and only if $\phi$ is a polynomial of degree at most $k-1$. The ``if'' direction of this assertion may be established by induction on the degree, since if $\phi$ is a polynomial of degree $d$ then, for fixed $h$, $\Delta_h \phi$ is a polynomial of degree $d-1$. Then ``only if'' direction can then be established by taking $h_1 = \dots = h_k = 1$ in \eqref{add-der-cond}, which tells us that $\phi(x + k)$ is uniquely determined as a function of $\phi(x), \phi(x+1), \dots, \phi(x+k - 1)$. Therefore $\phi$ is uniquely determined by its values at $0,1,\dots, k-1$, and hence coincides with the unique polynomial of degree at most $k-1$ which agrees with it at those points. 

The stability question, that is to say the characteristation of those $f$ for which $v(f) \geq 1- o(1)$, is already interesting. It turns out that $f$ must be closely approximated by a polynomial phase $e^{2\pi i \phi(x)}$. A precise statement and proof of this result may be found in \cite[Theorem 1.2]{eisner-tao}. The argument there is analogous to an earlier argument \cite{AKKLR} in a finite field setting, which has applications to property testing in theoretical computer science.

As with Young's inequality, however, our main focus here will be on the robustness question: for which $f$ do we have $v(f) \geq \frac{1}{K}$? This is known as the \emph{inverse question for the Gowers norms}. When $k = 2$, all such $f$ are at least somewhat related to exponentials of linear phases (the solutions to the equality question $v(f) = 1$).

\begin{theorem}\label{u2-inverse}
Suppose that $f : [N] \rightarrow \C$ is a function with $|f(x)| \leq 1$ for all $x$, and that $\Vert f \Vert_{U^2[N]} \geq \frac{1}{K}$. Then there is some $\theta \in \R/\Z$ such that \[ \frac{1}{N} |\sum_{x \in N} f(x) e^{-2\pi i \theta x}| \geq \frac{1}{K^2}.\] 
\end{theorem}

The proof of this is an exercise in Fourier analysis, given in detail in \cite[Proposition 8.2]{green-icm2006}. When $k \geq 3$, however, the situation is different. There are examples of functions $f : [N] \rightarrow \C$ with $|f(x)| \leq 1$ for all $x$, $\Vert f \Vert_{U^3[N]} \geq \frac{1}{K}$, but for which \begin{equation}\label{eq9} \frac{1}{N} | \sum_{x \in N} f(x) e^{-2\pi i \phi(x)} | \ll N^{-c}\end{equation} for all quadratic phases $\phi : \Z \rightarrow \R/\Z$. It is actually rather easy to give an example of such a function, though considerably less easy to prove rigorously that it \emph{is} an example: take $f(x) = e^{2\pi i \alpha x \{ \beta x\}}$, with $\alpha, \beta \in \R$ sufficiently irrational numbers such as $\alpha = \sqrt{2}$ and $\beta = \sqrt{3}$. Here $\{ t \}$ denotes fractional part. The ``reason'' this function $f$ has large $U^3[N]$-norm is that the phase $\phi(x) = \alpha x \{ \beta x\}$, whilst it does not satisfy the derivative condition \eqref{add-der-cond} exactly, does satisfy this condition for a positive proportion of $x, h_1, h_2, h_3$: in fact whenever $\{\beta x\}, \{ \beta h_1\}, \{\beta h_2\}, \{\beta h_3 \} \in [-\frac{1}{10}, \frac{1}{10}]$. Establishing \eqref{eq9} rigorously is quite tricky.

A more natural way to construct such functions is as \emph{nilsequences}. These objects should be thought of as ``higher-order characters'' generalising the linear exponentials $\Phi(n) = e^{2\pi i \theta  n}$. To explain this generalisation we write $\Phi$ in the form
\begin{equation}\label{chi-def} \Phi(n) = F(p(n))\end{equation}
where
\begin{itemize}
\item $p(n) = T^n 0$, where $T : \R \rightarrow \R$ is the translation map $Tx = x + \theta$;
\item $F(x) = e^{2\pi i x}$. Note that this function is $\Z$-periodic. 
\end{itemize}
A nilsequence corresponds to a generalisation of this in which $\R$ is replaced by a simply-connected nilpotent Lie group $G$ and $\Z$ is replaced by a lattice $\Gamma \subset G$. With this setup, a nilsequence is of the form \eqref{chi-def} with
\begin{itemize}
\item $p(n) = T^n \id_G$, where $T : G \rightarrow G$ is a \emph{nilrotation}, that is to say a map of the form $Tx = xg$ for some $g \in G$;
\item $F : G \rightarrow \C$ is smooth and $\Gamma$-\emph{automorphic}, which means that $F(\gamma x) = F(x)$ for all $x \in G$ and $\gamma \in \Gamma$. 
\end{itemize}
For example, we could take $G$ to be the Heisenberg group $\heis{\R}{\R}{\R}$ and $\Gamma$ to be the lattice $\heis{\Z}{\Z}{\Z}$. In fact for various reasons one usually considers a generalisation of this in which $p(n)$ is a \emph{polynomial sequence} on the group $G$. We will not discuss this important issue here, save to remark that it leads to essentially the same concept in the end due to a lifting argument of Furstenberg \cite[p. 31]{furst-book} (see also \cite[Appendix C]{gtz}). We say that $\Phi$ is an $s$-step nilsequence if the underlying nilpotent group $G$ has nilpotency class $s$, that is to say if the lower central series of $G$ is 
\[ G = G_1 \supset G_2 \supset \dots \supset G_s \supset G_{s+1} = \{\id_G\},\] with $G_s$ nontrivial. For the Heisenberg group we have $s = 2$.

To give a specific example in the Heisenberg case, we need to specify $g$ and an automorphic function $F$. The element $g$ can of course be specified just by  choosing a matrix. Nontrivial automorphic functions can be defined by hand (define $F$ to be a smooth bump function supported on the interior a fundamental domain for $\Gamma \backslash G$ and extend by automorphy). In the Heisenberg case there is a construction, pointed out in \cite{host-kra-nilclassification} for example, using the \emph{Jacobi $\theta$-function} $\theta(u,z) := \sum_n e^{\pi i z n^2 + 2\pi i nu}$ by defining $F( \heis{x}{y}{z}) = e(z) e^{-\pi x^2} \theta(y + ix)$. 

Functions of the form $e^{2\pi i \phi(x)}$ with $\phi$ quadratic (that is to say, the functions for which $\Vert f \Vert_{U^3[N]} = 1$) ``morally'' arise as nilsequences on the Heisenberg group by taking $g = \heis{1}{\alpha}{0}$ and $F(x,y,z) = e^{2\pi i (z - y \lfloor x\rfloor)}$. Indeed it may be checked by a computation that in this case we have $\Phi(n) = e^{2\pi i \phi(n)}$ with $\phi(n) = -\frac{1}{2}\alpha n (n+1)$. The slight technical issue here is that, although $F$ is automorphic (as may be confirmed by a computation) it is only piecewise smooth.

The following turns out to be true.

\begin{theorem}\label{robust-gowers}
Suppose that $\Phi(n)$ is an $s$-step nilsequence with\footnote{Here $\Vert \Phi \Vert_2^2 = \frac{1}{N} \sum_{n \leq N} |\Phi(n)|^2$. Some condition is needed to ensure that we do not have $\Phi(n) = 0$ identically.} $\Vert \Phi \Vert_2 = 1$. Then $\Vert \Phi \Vert_{U^{s+1}[N]} \geq \frac{1}{K}$, where $K$ is bounded above in terms of $s$ and the ``complexity'' of $\Phi$. 
\end{theorem}

Giving a proper definition of the complexity is a rather tedious matter; it must take account of various parameters associated with $G, \Gamma$ and the smoothness of the automorphic function $F$. The appendices of \cite{gt-nilratner} go into considerable further detail. 

The key to the proof is the observation that the multiplicative derivative $\Delta_h \Phi(n) = \Phi(n) \overline{\Phi(n+h)}$ is an $(s-1)$-step nilsequence, which allows us to proceed inductively. In fact, this is not quite true, but it is true if the automorphic function $F$ has the additional transformation property 
\begin{equation}\label{central-invariance} F(g_s x) = \xi(g_s) F(x)\end{equation} for every $g_s$ in $G_s$, the last nontrivial subgroup in the lower central series of $G$, for some character $\xi : G_s \rightarrow \C$ invariant under $\Gamma$. One may reduce to this case by a Fourier expansion on cosets of $G_s$. It is in effecting this Fourier expansion that the complexity of $\Phi$, and in particular the smoothness properties of $F$, comes into play. Suppose now that we do have the transformation property \eqref{central-invariance}. For fixed $h$ we have 
\[ \Delta_h \Phi(n) = F(T^n \id_G) \overline{F(T^{n+h} \id_G)} = \tilde F(T^n \id_{G})\] where  $\tilde F(x) = F(x) \overline{F(T^h x)}$. The function $\tilde F$ is easily seen to be $\Gamma$-automorphic, and moreover it is invariant under $G_s$:
\[ \tilde F(g_s x) = F(g_s x) \overline{F(T^n g_s x)} = \xi(g_s) F(x) \overline{\xi(g_s)} \overline{F(T^n x)} = F(x) \overline{F(T^n x)} = \tilde F(x).\] Here we used the fact that $G_s$ is central in $G$ to commute $T$ with multiplication by $g_s$. As a consequence, $\tilde F$ decends to an automorphic function on $G/G_s$, a nilpotent Lie group of class $s-1$. (Unfortunately the preceding discussion was actually quite a serious oversimplification, as in the definition of Gowers norm $h$ is not fixed but can vary over $[-N, N]$. With the argument just described, various smoothness norms of $\tilde F$ depend heavily on $h$ and to get around this a more complicated construction is required. Such a construction is given in \cite[Section 7]{gt-nilratner}.)

We have seen that functions supplying equality in \eqref{gowers-triv}, the inequality $v(f) \leq 1$, are polynomial phases. Theorem \ref{robust-gowers} states that nilsequences of step $k-1$ and suitably bounded complexity are solutions to the corresponding robustness problem $v(f) \geq \frac{1}{K}$, and so we think of them as ``approximate polynomials'' (or, more accurately, approximate polynomial phases). The discussion of the previous paragraph, where we saw that the multiplicative derivatives of nilsequences (with an additional invariance property) are nilsequences of lower step, adds further weight to this philosophy. 

It is very far from true that every solution to the robustness problem is a nilsequence. Indeed\footnote{Passing from $f$ to $f + \eps$ may destroy the property $|f(x)| \leq 1$, but we ignore this for the sake of illustration.} if $f : [N] \rightarrow \C$ and if $\eps : [N] \rightarrow \{-1,1\}$ is a random $\pm 1$-valued function then almost surely $v(f) \approx v(f + \eps)$. However, it is true that every solution is somewhat related to a nilsequence. 

\begin{theorem}\label{gis-thm}
Suppose that $f : [N] \rightarrow \C$ is a function with $|f(n)| \leq 1$ and that $\Vert f \Vert_{U^k[N]} \geq \frac{1}{K}$. Then there is a $(k-1)$-step nilsequence $\Phi(n)$ with $|\Phi(n)| \leq 1$ for all $n$ such that 
\[ \frac{1}{N} |\sum_{n \in N} f(n) \overline{\Phi(n)}| \geq \frac{1}{K'},\]
where $K'$ and the complexity of $\Phi$ are bounded in terms of $K$ and $k$ only.
\end{theorem}
Note that Theorem \ref{gis-thm} is a generalisation of Theorem \ref{u2-inverse}, which was essentially the case $k = 2$. This result is due to Tao, Ziegler and the author \cite{gtz} and is known as the \emph{Inverse Theorem for the Gowers norms}. A weaker ``local'' version of it was obtained by Gowers (in \cite{gowers-4aps} for $k = 3$, and in \cite{gowers-longaps} for general $k$). The case $k = 3$ was established by Tao and the author \cite{gt-u3inverse}, and the case $k = 4$ by Tao, Ziegler and the author \cite{gtz-u4}. It should most certainly be mentioned that the relevance of nilpotent Lie groups in this general arena first became apparent in the context of ergodic theory in works of Conze, Lesigne, Furstenberg and Weiss \cite{conze-lesigne,furst-survey-old,furst-survey,furst-weiss}.  A result which may be thought of as an ``ergodic analogue'' of Theorem \ref{gis-thm} was obtained by Host and Kra \cite{host-kra} (see also independent work of Ziegler \cite{ziegler}). The notion of nilsequence itself, as well as the word, was introduced by Bergelson, Host and Kra \cite{bhk}. See \cite{kra-icm} for a nice introduction to these connections.

The conjecture to which Theorem \ref{gis-thm} is a solution, together with potential applications of it to prime numbers, was formulated by Tao and the author \cite{gt-linearprimes} about four years before it was finally proved. The proof is unfortunately inordinately long and complicated. For a summary in about 20 pages, see \cite{gtz-announce}. An alternative approach has been developed by Szegedy \cite{szegedy} and Camarena-Szegedy \cite{camarena-szegedy}, based in part on the work of Host and Kra, but these papers are not an easy read either. 

The converse of Theorem \ref{gis-thm} is also true, with appropriate changes to the constants. The proof is relatively straightforward and goes along very similar lines to the sketch of the proof of Theorem \ref{robust-gowers} we gave above: see \cite[Appendix G]{gtz-u4} for further details.

Although we do not plan to discuss it much here, there has also been a good deal of work on finite field analogues of Theorem \ref{gis-thm}, which have applications in theoretical computer science. In addition to work by various subsets of the authors named above, we note that Samorodnitsky \cite{samorod} established the case $k = 3$ of Theorem \ref{gis-thm} in the particularly interesting setting where $[N]$ is replaced by a vector space in characteristic 2. 

\begin{theorem}\label{samorod}
Suppose that $V$ is a finite-dimensional vector space in characteristic 2 and that $f : V \rightarrow \C$ is a function with $|f(x)| \leq 1$ for all $x \in V$. Suppose that $\Vert f \Vert_{U^3(V)} \geq \frac{1}{K}$. Then there is a function $\Phi : V \rightarrow \C$ of the form $\Phi(x) = (-1)^{\psi(x)}$, where $\psi : V \rightarrow \F_2$ is a quadratic form, such that 
\[ \frac{1}{|V|} |\sum_{x \in V} f(x) \Phi(x)| \geq \frac{1}{K'}.\]
Here $K'$ is bounded in terms of $K$ only.
\end{theorem}
The definition of the $U^3(V)$-norm is entirely analogous to that of the $U^3[N]$ norm, except that the average $\E$ is now simply taken over all $x, h_1, h_2, h_3 \in V$. Samorodnitsky obtained a bound of the form $K' \sim e^{K^C}$, but by combining his methods with the work of Sanders \cite{sanders-bogolyubov} one could improve this to $K' \sim e^{(\log K)^C}$. The similarity of these bounds to those stated in conjunction with Theorem \ref{ruz-fin-thm} is no coincidence. Indeed a close relationship between the structure theory of approximate subgroups of $\F_2^{\Z}$ and Theorem \ref{samorod} was discovered by Tao and the author \cite{gt-cambphilsoc} and independently by Lovett \cite{lovett}. In particular, it is known that the Polynomial Freiman-Ruzsa conjecture for finite fields, which is equivalent to Conjecture \ref{conj2.9}, is also equivalent to having a bound of shape $K' \ll K^C$ in Theorem \ref{samorod}. 

A similar equivalence between bounds in Theorem \ref{freiman-ruzsa-theorem} and the case $k = 3$ of Theorem \ref{gis-thm} was developed in \cite{gt-cambphilsoc}: in other words the theories of approximate subgroups of $\Z$ and of approximate quadratic polynomials are in a sense the same. I have often informally advanced the speculation that looking for a more effective proof of Theorem \ref{gis-thm} may be a way of attacking the Polynomial Freiman-Ruzsa Conjecture, though without any convincing ideas about how this might be achieved.

\subsection{Applications} The theory of Gowers norms as described in the last section was for the most part developed to study arithmetic progressions. Gowers himself was interested in Szemer\'edi's theorem, and Tao and the author were subsequently concerned with arithmetic progressions of primes. In \cite{gt-linearprimes} it was observed that the theory applied to a fairly wide class of ``linear'' problems, including questions about linear configurations of primes. Since this theory was discussed\footnote{Naturally, however, this account is quite out of date and in particular predates the general case of Theorem \ref{gis-thm}.} in my 2006 ICM lecture \cite{green-icm2006} and is described in the article of Ziegler in these \emph{Proceedings}, we restrict ourselves to a very brief account.

The connection of the Gowers norms to linear configurations comes from results called \emph{generalised von Neumann inequalities}, which have the form
\begin{equation}\label{gvn} | T(f_1,\dots, f_t)| \ll \inf_{i = 1,\dots, t}\Vert f_i \Vert_{U^{s+1}[N]}.\end{equation}
Here, $f_1, \dots, f_t : [N] \rightarrow [-1,1]$ are functions and 
\[ T(f_1,\dots, f_t) = \E_{(n_1,\dots, n_d) \in S} f_1(\psi_1(n_1,\dots, n_d)) \dots f_t(\psi_t(n_1,\dots, n_d)),\] where the $\psi_i : \Z^d \rightarrow \Z$ are affine-linear forms and $S$ is a ``nice'' set (for example a convex set). For any system of forms $\psi_1,\dots, \psi_d$ which is not degenerate in a certain sense, there is a value of $s$ for which \eqref{gvn} holds. For example, if $d = 2$, $t = 3$ and $\psi_1(n_1, n_2) = n_1$, $\psi_2(n_1, n_2) = n_1 + n_2$, $\psi_3(n_1, n_2) = n_1 + 2n_2$ (3-term arithmetic progressions) then we may take $k = 2$, whilst if $d = 2$, $t = 4$ and $\psi_1(n_1, n_2) = n_1$, $\psi_2(n_1, n_2) = n_1 + n_2$, $\psi_3(n_1, n_2) = n_1 + 2n_2$, $\psi_4(n_1, n_2) = n_1 + 3n_2$ (4-term arithmetic progressions) then we may take $k = 3$. The degenerate confirgurations are those in which some two of the $\psi_i$ have equal homogeneous part, up to scalar equivalence: thus for example we cannot take $\psi_1(n_1) = n_1$ and $\psi_2(n_1) = n_1 + 2$.  The proof of any generalised von Neumann inequality is conceptually quite easy, involving only several applications of the Cauchy-Schwarz inequality, but notationally a little unpleasant. A general form of \eqref{gvn} was obtained in \cite[Appendix D]{gt-linearprimes}. Furthermore the inequality was established there under a weaker condition on the $f_i$ than boundedness, namely that $|f_i(x)| \leq \nu(x)$ pointwise for some ``pseudorandom measure'' $\nu$. This is crucial for applications to the primes. 

Ignoring a few technicalities, the manner in which \eqref{gvn} is applied to the primes is as follows. For technical convenience the primes are weighted using the von Mangoldt function $\Lambda$, defined by $\Lambda(n) = \log p$ if $n = p^k$ is a prime power and $\Lambda(n) = 0$ otherwise. We are interested in $T(\Lambda, \dots, \Lambda)$, which counts how often the linear forms $\psi_1(n_1,\dots, n_d), \dots, \psi_t(n_1,\dots, n_d)$ all take prime values as $(n_1,\dots, n_d)$ ranges over a set $S$. 
To estimate this we split $\Lambda$ in a certain manner as 
\begin{equation}\label{decomp} \Lambda = \Lambda^{\sharp} + \Lambda^{\flat},\end{equation} where $\Lambda^{\sharp}$ is ``structured'' and $\Lambda^{\flat}$ is ``unstructured''. Since $T$ is multilinear, we may split $T(\Lambda, \dots, \Lambda)$ as a sum of $T(\Lambda^{\sharp}, \dots, \Lambda^{\sharp})$ plus $2^t - 1$ other terms, each of which involves at least one copy of $\Lambda^{\flat}$. The first term provides the main term in the asymptotic formula for $T(\Lambda, \dots, \Lambda)$, and the aim is then to show that the other $2^t - 1$ terms are all small. By \eqref{gvn}, this may be accomplished if it can be shown that 
\[ \Vert \Lambda^{\flat} \Vert_{U^{s+1}[N]} = o(1).\]
By the inverse theorem for the Gowers norms, Theorem \ref{gis-thm} (in the contrapositive), it is enough to establish that
\begin{equation}\label{task} \frac{1}{N} |\sum_{n \leq N} \Lambda^{\flat}(n) \overline{\Phi(n)}| = o(1)\end{equation} for every $s$-step nilsequence $\Phi(n)$ of bounded complexity. At least, this would be so were it not for the restriction $|f(x)| \leq 1$ in Theorem \ref{gis-thm}: a large part of \cite{gt-linearprimes} is devoted to removing this restriction, showing that Theorem \ref{gis-thm} implies a more general version of itself in which we only assume that $|f(x)| \leq \nu(x)$ for some pseudorandom measure $\nu$. 

The actual decomposition \eqref{decomp} we choose is based on the formula $\Lambda(n) = \sum_{d | n} \mu(d) \log(n/d)$, where $\mu$ is the M\"obius function. It transpires that the task of establishing \eqref{task} may be further reduced to establishing that 
\begin{equation}\label{mobius-nil} \frac{1}{N} |\sum_{n \leq N} \mu(n) \overline{\Phi(n)}| \ll_A \log^{-A} N\end{equation} for every $A > 0$. This statement was formerly known as the ``M\"obius and nilsequences conjecture'', but it is now a theorem of Tao and the author \cite{gt-nilmobius}. Although the paper \cite{gt-nilmobius} is relatively short, it depends crucially on the much longer paper \cite{gt-nilratner}, in which various properties of nilsequences are established, in particular with regard to the distribution of finite orbit segments $(T^n\id_G)_{n \leq N}$ in $\Gamma \setminus G$. This work, like other material in this section, was motivated by earlier developments in the ergodic theory community, in particular work of Leon Green \cite{leon-green} and papers of Leibman of both an algebraic \cite{leibman1} and an ergodic-theoretic \cite{leibman2} nature.

\subsection{Open questions} For me the key open question is to find the ``right'' proof of the inverse conjecture for the Gowers norms. At the moment the proofs are unsatisfactory on a conceptual level (the notion of a nilsequence is extremely natural, so it would be disappointing if it genuinely required 100+ pages to explain its role in Theorem \ref{gis-thm}). Furthermore, these proofs provide rather poor bounds for the complexity of the nilsequence $\Phi$, particularly when $k \geq 4$ (in fact for $k \geq 5$ the proofs provide no explicit bounds at all due to the use of ultrafilter arguments, though once again an explicit bound could in principle be extracted via quantifier elimination). As noted above it would be particularly interesting, in view of the link to approximate subgroups of $\Z$, to find a new approach to the inverse theorem when $k = 3$.

A more specific question is whether there is some smaller ``natural'' class of nilsequences. The space $C^{\infty}(\Gamma \setminus G)$ of automorphic functions is extremely large, but we know for example that in the case $\Gamma = \Z$, $G = \R$ the exponentials $e^{2\pi i x}$ have a special role. Eigenfunctions of Laplacians are one natural avenue of enquiry. Furthermore the space of all simply-connected nilpotent Lie groups $G$ together with lattices $\Gamma$ is also extremely large and complicated, and it may be natural to focus on some subclass (for example free nilpotent Lie groups).

\section{Other directions}

To conclude this article I want to mention a personal selection of a few other inequalities where the equality, stability and robustness questions may hide interesting algebraic or somewhat algebraic structure. In some cases there is at least a tenuous connection to the main sections of the article, and in others less so.

\subsection{Inverse questions for the large sieve}  Let $\mathscr{A}$ be a set of natural numbers with the property that $|\mathscr{A} \pmod{p}| \leq \frac{1}{2}(p+1)$ for all sufficiently large primes $p$. The large sieve guarantees that $|\mathscr{A} \cap [N]| \ll N^{1/2}$ for all $N$. This is sharp up to a multiplicative constant, as is shown by taking $\mathscr{A}$ to be the set of squares (or the set of integer values of an arbitrary quadratic with rational coefficients). 

It may well be the case that a very strong robustness assertion holds: if there is some $K$ such that $|\mathscr{A} \cap [N]| \geq \frac{1}{K} N^{1/2}$ for all sufficiently large $N$ then $\mathscr{A}$ is contained, up to a finite set, in the set of values of a rational quadratic. See \cite{green-harper} for evidence in this direction. This type of question was first raised by Helfgott and Venkatesh \cite{hv}; see also \cite{walsh1, walsh2}.

\subsection{Point-line configurations} Let $\mathscr{P} \subset \R^2$ be a set of $n$ points, no four on a line\footnote{This condition is included here for simplicity, but can probably be relaxed.}. Write $T(\mathscr{P})$ for the number of pairs $(x,y) \in \mathscr{P}$ of distinct points for which there is a third distinct point $z \in \mathscr{P}$ on the line $\overline{xy}$. Trivially, $T(\mathscr{P}) \leq n(n-1)$. Less obviously, equality cannot occur: this follows from a famous result known as the Sylvester--Gallai theorem. 

Almost-equality can occur: we can obtain $T(\mathscr{P}) = n^2 - O(n)$ by taking $\mathscr{P}$ to be a suitable set of points on a suitable cubic curve (for example a coset of a subgroup on an elliptic curve, although there are singular examples too). This was noted by Sylvester in the 1860s \cite{sylvester}. Conversely, it was recently shown by Tao and the author \cite{green-tao-pointslines} that there is a strong converse to this statement. 

It would be very interesting to have an understanding of those $\mathscr{P}$ for which $T(\mathscr{P}) = n^2(1 - o(1))$ (the stability question) or, more ambitiously, $T(\mathscr{P}) \geq \frac{n^2}{K}$ (the robustness question). The paper \cite{green-tao-pointslines} only covers the extreme end of the stability region. It is possible that cubic structure is responsible for all such $\mathscr{P}$. There are links here to the theory of approximate groups: for example, finite approximate subgroups of elliptic curve groups are a source of examples of such sets $\mathscr{P}$. 

An interesting nontrivial result in higher dimensions is \cite{wigderson-et-al}, motivated by applications in theoretical computer science.

\subsection{The Littlewood Problem} Suppose that $A \subset \Z$ is a set of $n$ integers. Then it was established 30 years ago by Konyagin \cite{konyagin} and McGehee-Pigno-Smith \cite{mps}, answering a question of Littlewood \cite{littlewood}, that 
\[ \int^1_0 |\sum_{a \in A} e^{2\pi i \theta a}| d\theta \gg \log n.\]
Earlier results had been obtained by Paul Cohen and others. This is sharp up to the constant, as is shown by taking $A$ to be an arithmetic progression of length $n$. (In fact, this example may also provide the sharp constant, a conjecture known as the \emph{Strong Littlewood Conjecture}.) Very little is known about the robustness question, that is to say about the structure of those $A$ for which
\[ \int^1_0 |\sum_{a \in A} e^{2\pi i \theta a}| d\theta \leq K \log n.\] It is possible that such $A$ are very close to being unions of a few arithmetic progressions. If so, this would have applications to questions in combinatorial number theory about sum-free sets due to a connection established by Bourgain \cite{bourgain-sum-free}. For some partial results and a further discussion, see \cite{petridis}.

\newcommand\PG{\operatorname{PG}}
\subsection{No-three-in-a-line} Let $p$ be an odd prime, and suppose that $A \subset \PG(2,p)$ is a set containing no three distinct points in a line\footnote{Such sets are called ``arcs'' in the literature, which is extremely extensive.}. (Here, $\PG(2,p)$ is the 2-dimensional projective space over $\F_p$, thus $|\PG(2,p)| = p^2 + p+ 1$). It is very easy to see that $|A| \leq p+2$ and an exercise to show that $|A| \leq p+1$. Equality occurs when $A$ is a conic.  Remarkably, a celebrated result of Segre \cite{segre} shows that in fact equality occurs \emph{only} when $A$ is a conic. 

The stability question was resolved by Voloch \cite{voloch}, building upon remarkable work of Segre. Voloch shows that any $A$ with no three-in-a-line and $|A| \geq \frac{44}{45} p$ is contained in a conic. This argument is quite deep, depending on an application of the polynomial method \cite{tao-polynomial-method} as well as bounds of St\"ohr and Voloch \cite{stohr-voloch} about counting points on high degree curves. 

The robustness question, that is to say the classification of those $A$ with $|A| \geq \frac{1}{K}p$, is very interesting. There are examples coming from cubic curves, such as $A = \{(x : x^3 : 1) : 0 <x < p/3\}$. So far as I am aware there is no example in the literature to contradict the possibility that all sets $A \subset \PG(2,p)$ with no-three-on-a-line and $|A| \geq \frac{1}{K}p$ have all but $o(p)$ of their points lying on a curve of degree at most $3$. So far as I am aware no-one has explicitly conjectured this either, so perhaps I shall take this opportunity to do so. 

There is a superficial link to a notorious problem of Dudeney \cite{dudeney} about whether there is a set $A$ of $2N$ points on the grid $[N] \times [N]$ with no three in a line. There are many fewer colinear triples in $[N] \times [N]$ than in $\Z/p\Z \times \Z/p\Z$ for $p \sim N$, however, so the study of sets $A$ such as this is likely to be even harder than the problem discussed above. Nonetheless, the best-known examples (with $|A| \sim 3N/2$, see \cite{hjsv}) are given by very algebraic constructions. It seems likely that the answer to Dudeney's question is negative.

\subsection{Sidon sets} Suppose that $A \subset [N]$ is a set with the property that all pairwise sums $x + y$ with $x, y \in A$ are distinct, apart from the obvious coincidences $x + y = y + x$. Such a set $A$ is called a Sidon set. It is very easy to see that $|A| \ll \sqrt{N}$, and with more care (an argument of Erd\H{o}s and Tur\'an) one may show that $|A| \leq (1 + o(1)) \sqrt{N}$. There are examples of Sidon sets $A$ with $|A| = (1 - o(1))\sqrt{N}$, all constructed in a highly algebraic maner using finite fields. There are different variants due to Bose, Ruzsa and Singer. It is possible that the stability question (that is, the classification of those Sidon set $A$ with $|A| = (1 - o(1))\sqrt{N}$) has a satisfactory answer, but there is no obvious guess, based on the known examples, as to what it might be. The robustness question, that is to say the classification of those $A$ with $|A| \geq \frac{1}{K}\sqrt{N}$, is of course even more difficult. A discussion of it was had on the blog of Tim Gowers \cite{gowers-blog}. In commenting on that discussion, Terence Tao raised the possibility that an answer to this question could lead to progress on a famous and old problem of Erd\H{o}s, namely to determine if there is an additive basis $\mathscr{A}$ of the natural numbers of order $2$ (i.e. $\mathscr{A}+\mathscr{A}=\mathbb{N}$) with an absolute bound on the number of representations of $x$ as a sum of two elements of $\mathscr{A}$.

\section{Acknowledgements}

I would like to thank Sean Eberhard, Bryna Kra, Freddie Manners, Peter Sarnak and Terence Tao for comments on a draft of this article. I wish to thank the last of these for our extensive collaboration over the last decade, which has so far led to 30 joint papers. I also thank my other coauthors, three of whom are speaking at this congress. It is only because of these mathematicians that I have the opportunity to present these topics at the 2014 ICM.

\end{document}